\documentclass[12pt]{arxivclassfile}

\def\R{{\mathbb R}}

\def\dee{{\rm d}}
\def\e{{\rm e}}

\newcommand{\p}{\partial}
\newcommand{\ft}{\tilde{f}}

\newcommand{\A}{\alpha}
\newcommand{\B}{\beta}

\newcommand{\be}{\begin{equation}}
\newcommand{\bi}{\begin{itemize}}
\newcommand{\ei}{\end{itemize}}
\newcommand{\ee}{\end{equation}}
\newcommand{\ba}{\begin{eqnarray}}
\newcommand{\bs}{\begin{eqnarray*}}
\newcommand{\ea}{\end{eqnarray}}
\newcommand{\es}{\end{eqnarray*}}

\newcommand{\Proof}{\noindent {\it Proof. }}
\newcommand{\EndProof}{\hfill $\Box$\newline }

\usepackage{mathscinet}
\usepackage{url}
\usepackage{amssymb}
\usepackage{amsmath}
\usepackage{amsthm}
\usepackage{enumerate}
\newtheorem{theorem}{Theorem}[section]

\newtheorem{lemma}{Lemma}[section]
\newtheorem{proposition}{Proposition}[section]

\newtheorem{definition}{Definition}[section]

\begin{document}

\begin{frontmatter}


\author{R.\ Laister\corref{thing}\fnref{uwe}}
\ead{Robert.Laister@uwe.ac.uk}
\fntext[uwe]{Permanent address: Dept. of Engineering Design and Mathematics, University of the West of England, Bristol BS16 1QY, UK.}
\author{J.C.\ Robinson}
\ead{J.C.Robinson@warwick.ac.uk}
\author{M.\ Sier{\.z}\polhk{e}ga}
\ead{M.L.Sierzega@warwick.ac.uk}
\cortext[thing]{Corresponding author}

\title{Non-existence of local solutions for semilinear heat equations of Osgood type}



\address{Mathematics Institute, Zeeman Building, University of Warwick,\\ Coventry CV4 7AL, UK}

\begin{abstract}
We establish non-existence results for the Cauchy problem of some  semilinear heat equations
with non-negative initial data  and locally Lipschitz, nonnegative source term $f$. Global (in time) solutions of the scalar ODE $\dot v=f(v)$ exist for $v(0)>0$ if and only if the Osgood-type condition $\int_{1}^{\infty}\frac{\dee s}{f(s)} =\infty$ holds; by comparison this ensures the existence of global  classical solutions of $u_t=\Delta u+f(u)$ for bounded initial data $u_0\in L^{\infty}(\R^n)$. It is natural to ask whether the Osgood condition is sufficient to ensure that the problem still admits global solutions if the initial data is in $L^q(\R^n)$ for some $1\le q<\infty$. Here we answer this question in the negative, and in fact show that there are initial conditions for which there exists no local solution in $L^1_{\rm loc}(\R^n)$ for $t>0$.

\end{abstract}


\end{frontmatter}


\section{Introduction}
\label{Intro}

We consider the Cauchy problem for the  semilinear parabolic equation
\be
 u_t = \Delta u +  f(u),\qquad u_0\in L^q(\R^n )\label{eq:PDE}
 \ee
under the following hypotheses for the initial data $u_0$ and source term $f$:
\bi
\item[(H1)] $ u_0\in L^q(\R^n )$  and $u_0\ge 0$;
\item[(H2)] $f:[0,\infty )\to [0,\infty )$ is locally Lipschitz continuous, non-decreasing, $f(0)=0$,  $f>0$ on $(0,\infty )$;
\item[(H3)] $f$ satisfies the Osgood-type condition
\be
\int_{1}^{\infty}\frac{\dee s}{f(s)} =\infty .\label{eq:OC}
\ee
\ei
When (H2) holds, it is well-known that solutions of the ODE
\be
\dot{v}=f(v), \qquad v(0)>0,\label{eq:ODE}
\ee
exist globally in time if and only if (H3) holds. For   non-negative $ u_0\in L^\infty (\R^n )$ and $f$ satisfying (H2), it is also well-known  that a global (in time) solution of (\ref{eq:PDE}) exists whenever  (H3) holds. In this case a standard local existence argument (e.g.\ \cite{Henry,Pazy}) yields a solution that is bounded for a short time. By means of comparison with the  solution of the ODE (\ref{eq:ODE}) with initial condition $v(0)=\|u_0\|_{\infty}$,
this solution can then be extended for all times. A natural question to ask therefore, is whether (H3) still guarantees the global existence of solutions of (\ref{eq:PDE})  with  singular, {\em unbounded} initial data (see \cite{RS1}).

It is surprising  that this seemingly obvious question has remained open until now. The present work will resolve this question in the negative; indeed, we show that for a family of $f$ satisfying (H1--H3) there is initial data in $L^q(\R^n)$ for which there is no local solution of (\ref{eq:PDE}) in the following very general sense.

\begin{definition}{\cite[p.\ 78]{SQ}}
Given $f$ non-negative and $u_0\ge0$ we say that $u$ is a {\rm{local integral solution}} of (\ref{eq:PDE}) on $[0,T)$ if $u:\R^n\times[0,T)\to[0,\infty]$ is measurable, finite almost everywhere, and
\be
u(t)=S(t)u_0+\int_0^t S(t-s)f(u(s))\, \dee s\label{eq:mild}
\ee
holds almost everywhere in $\R^n\times[0,T)$, where
\bs
(S(t)u_0)(x)= (4\pi t)^{-n/2}\int_{\R^n}\e^{-|x-y|^2/4t}u_0(y)\, \dee y
\es
is a classical solution of the linear heat equation
\bs
w_t = \Delta w ,\qquad w(0)=u_0.
\es\label{def:soln}
\end{definition}
We remark that non-existence of a solution in the sense of Definition~\ref{def:soln} implies the non-existence of mild solutions and of classical solutions \cite[p.\ 77--78]{SQ}.

 In some of the literature the non-existence of a local solution is known as {\it instantaneous blow-up} (\cite{CZ,PV,Sam1}). To the best of our knowledge, instantaneous blow-up has never previously been shown for the semilinear heat equation (\ref{eq:PDE})  with data satisfying  (H1-H3).

Even when a global solution exists (for bounded initial data) and is asymptotically unbounded, we see that for any finite time interval the trajectory `explores' only a bounded subset of the range of $f$. Thus one of the main differences between considering bounded and unbounded data in the context of the Osgood-type condition is that singular data interacts with the whole range of the source term.

It should also be noted that the Osgood-type condition (\ref{eq:OC}) is not {\it necessary} for
local or global existence of solutions to (\ref{eq:PDE}), see \cite{BC1,BCMR,F1,NS}. Conversely, if $u_0\in L^q(\R^n)$ with $u_0\ge0$ then it is well known that the local existence of classical solutions of  (\ref{eq:PDE})  is guaranteed if $f$ satisfies the growth bound
\be
f(u)\le Cu^p\label{eq:locexist}
\ee
for some $p\in [1,1+2q/n)$, with $p=1+2q/n$ allowed if $q>1$; this follows using monotonicity methods, see \cite{RS2,WEISSLER1980,W} (for related results see \cite{BC1,GIGA1986,SQ}). An example of a non globally-Lipschitz function satisfying both (H2-H3) and (\ref{eq:locexist})
is $f(u)=(u+1)\ln{(u+1)}$, for which global existence of classical solutions is therefore assured.  We mention that Osgood-type conditions are  of interest in the study of well-posedness of aggregation equations \cite{BLR}.

Observe that the restriction imposed in (H3) by the Osgood-type condition is of an entirely different nature than the usual bounds on growth familiar from analysis of problems with power law nonlinearities like $f(u)=u^p$, $p>1$. A function $f$ may  display   arbitrarily fast growth locally without violating (H3). Furthermore, in many  semilinear problems one often finds a structural compatibility between the source term $f$ and the space of initial data. This compatibility plays a prominent role in calculations leading to well-posedness theory and blow-up analysis. An extensively studied  example  is the power law nonlinearity $f(u)=u^p$ and the space of initial conditions is modelled on Lebesgue spaces (\cite{BC1,GIGA1986,SQ,WEISSLER1979,WEISSLER1980}). Assumptions of this type naturally yield conditions expressed in terms of $p$ (growth of the source term), $q$ (integrability class of the initial data) and $n$, the dimension of the spatial domain. In principle, condition (\ref{eq:OC}) offers no such structural relation, ruling out a possibility of a unified approach to data taken in a  standard functional framework like $L^{p}$ or $W^{p,q}$.

Violation of (H3) has long been used for establishing sufficient conditions for blow-up in semilinear and quasilinear heat equations. In particular, when supplemented by a convexity assumption on $f$, it is the basis of a well-known eigenfunction argument \cite{Kaplan}. For a survey of such problems see \cite{GV} and for quasilinear heat equations see \cite{Sam1}. These works mainly address the problem of blow-up of solutions with bounded initial data.

The structure of the paper is as follows. In Section 2 we establish some preliminary estimates for the heat equation, and in Section 3 we construct a particular family of Osgood nonlinearities. Section 4 contains the main non-existence results for (\ref{eq:PDE}), and we conclude with some remarks and open problems in Section 5.


\section{Estimates for some solutions of the heat equation}
\label{main result}

For $r>0$,
$B_n(r)$ will denote the Euclidean ball in $\R^n$ of radius $r$, $S^{n-1}$ the unit sphere and $\omega_n$  the volume of the unit ball $B_n(1)$. For a particular choice of initial condition $u_0$ we find a set of spatial points at which $w(x,t)=(S(t)u_0)(x)$ is large.

\begin{proposition}\label{prop:ball}
Fix $\A\in (0,n)$ and $R>1$ and let $u_0\in L^1(\R^n) $ be the non-negative, radially symmetric function given by
$$
u_0(x)=|x|^{-\alpha}\chi_R:=\left\{\begin{array}{ll}
|x|^{-\A},& |x|\le R,\\
0,& |x|>R.
\end{array}\right.
$$
Let $w(t)=S(t)u_0$ and
\be
M=\min \{ w(\hat x,t)\ :\ \hat x\in S^{n-1},\ 0\le t\le 1\}.
\label{eq:M}\ee
If $\beta\in (0,1/2)$, then for any $\phi\ge M$
\be\label{wbig}
w(x,t)\ge \phi\quad\mbox{for}\quad |x|\le t^\beta,\qquad\mbox{for all}\quad 0< t\le (\phi/M)^{-1/\alpha\beta}.
\ee
\end{proposition}

\Proof
For any $0<t\le 1$ and $x=t^\B \hat{x}\in\p B_n(t^\B)$, $\hat{x}\in S^{n-1}$ we have
\begin{align*}
w(t^\B \hat{x},t) &= (4\pi t)^{-n/2}\int_{\R^n}\e^{-|t^\B \hat{x}-y|^2/4t}u_0(y)\, \dee y \\
&=(4\pi t)^{-n/2}\int_{\R^n}\e^{-t^{2\B} |\hat{x}-z|^2/4t}u_0(t^\B z)t^{n\B}\, \dee z \qquad\qquad (y=t^{\B}z)\\
& =   (4\pi )^{-n/2}t^{n\B -n/2-\A\B}\int_{B_n(Rt^{-\B})}\e^{-|\hat{x}-z|^2/4t^{1-2\B}}|z|^{-\A}\, \dee z \\
& \ge  (4\pi )^{-n/2} t^{n\B -n/2-\A\B}\int_{B_n({R})}\e^{-|\hat{x}-z|^2/4t^{1-2\B}}|z|^{-\A}\, \dee z\\
& =  t^{-\A\B}\left(4\pi t^{1-2\B}\right)^{-n/2} \int_{{\R^n}}\e^{-|\hat{x}-z|^2/4t^{1-2\B}}u_0(z)\, \dee z\\
& =  t^{-\A\B}w(\hat{x},t^{1-2\B}).
\end{align*}
Note that since $0<\B <1/2$, $0<t\le 1$ and $|\hat{x}|=1$, from (\ref{eq:M}) we have that
$$
w(\hat{x},t^{1-2\B})\ge M \qquad\mbox{for all}\qquad  \hat{x}\in S^{n-1},\ 0<t\le 1.
$$
Now observe that since $w$ is radially symmetric and decreasing in the radial variable,
$$
  w(x,t)\ge w(t^\beta \hat x,t)\ge Mt^{-\alpha\beta}\qquad\mbox{for all}\qquad |x|\le t^\beta,
  $$
  from which (\ref{wbig}) follows immediately.\EndProof

\section{A family of Osgood $f$}\label{sec:f}

In this section we construct a family of functions $f$ depending on a parameter $k>1$ that satisfy (H2-H3), suitable for   establishing our non-existence results for local solutions of (\ref{eq:PDE}).

For a given $k>1$, fix $\phi_0>2^{1/(k-1)}$ and define inductively the sequence $\phi_i$  via
$$
\phi_{i+1}={\phi_{i}^k}.
$$
It is easily verified that
\be
1<\phi_{i-1}<{\phi_{i}}/2\qquad \forall i\ge 1, \label{eq:phi}
\ee
and $\phi_i\to\infty$ as $i\to\infty$.

Now define $f:[0,\infty )\to [0,\infty )$ by
\ba
f(s)=\left\{
\begin{array}{ll}
(1-\phi_0^{1-k})s^k,& s\in J_0:=[0,{\phi_{0}}],\\
\phi_{i}-{\phi_{i-1}},& s\in I_i:=[\phi_{i-1},{\phi_{i}}/2],\ i\ge 1,\\
\ell_i(s),&s\in J_i:=(\phi_i/2,\phi_i),\ i\ge 1,
\end{array}\right.
\label{eq:f}
\ea
where $\ell_i$ interpolates linearly between the values of $f$ at $\phi_{i}/2$ and $\phi_i$. By construction $f$ is locally Lipschitz continuous, constant on the intervals $I_i$ (with value $\phi_{i}-{\phi_{i-1}}$) and linear (more precisely, affine) on the intervals $J_i$. Clearly $f$ satisfies (H2); it also satisfies the Osgood-type condition (H3), since
$$
\int_{1}^{\infty}\frac{\dee s}{f(s)}\ge \sum_i\int_{I_i}\frac{\dee s}{f(s)}= \sum_{i=1}^{\infty}\frac{\phi_{i}/2-{\phi_{i-1}}}{\phi_i-{\phi_{i-1}}}=\frac{1}{2}
\sum_{i=1}^{\infty}\left( 1-\frac{\phi_{i-1}}{\phi_i-\phi_{i-1}}\right)=\infty .
$$
We note also that our function $f(s)$ is bounded above by $2^ks^k$ for $s\ge 0$. This is clear for $s\in [0,\phi_0]$;  for $s\in I_i=[\phi_{i-1},\phi_i/2]$ we have
$$
f(s)=\phi_i-\phi_{i-1}\le \phi_{i-1}^k\le s^k;
$$
and for $s\in J_i=(\phi_i/2,\phi_i)$ we have
$$
f(s)\le\phi_{i+1}-\phi_i=\phi_i^k-\phi_i\le\phi_i^k\le 2^ks^k.
$$
We therefore obtain (see discussion around equation (\ref{eq:locexist})) global existence of solutions with initial data in $L^1(\R^n)$ provided that $k<1+2/n$; and global existence of solutions with initial data in $L^q(\R^n)$, $q>1$, provided that $k\le 1+2q/n$. However, we will see below that for any $k>q(1+2/n)$ there is initial data in $L^q(\R^n)$ for which there exists no local integral solution.

For comparison we also consider the function $\ft :[0,\infty )\to [0,\infty )$ defined by
$$
\ft (s)=\left\{
\begin{array}{ll}
0,& s\in [0,{\phi_{0}}]=J_0,\\
\phi_{i}-{\phi_{i-1}},& s\in (\phi_{i-1},{\phi_{i}}]=I_i\cup J_i,\ i\ge 1.
\end{array}\right.
$$
Then $\ft =f$ on $I_i$ and $f \ge \ft$ on $J_i$ and so  $f\ge\tilde f$ on $[0,\infty )$.


\section{Non-existence of local solutions}

Before we prove our non-existence theorem we give a quick proof of a weaker result that illustrates the essential idea: we show that if $k>1+2/n$ then there exists a $u_0\in L^1(\R^n)$ such that there is no local integral solution of (\ref{eq:PDE}) that remains in $L^1(\R^n)$. Note that in this case there is a sharp transition in the values of $k$ for which we have global existence ($k<1+2/n$) to non-existence for certain initial data. There is a similar sharp dichotomy in existence/non-existence of solutions  when $f(u)=u^p$ for initial data in $L^1(\R^n)$ \cite{BC1} and even for $u_0$ given by the Dirac-delta function \cite{BF}.

\begin{lemma}
  If $k>1+2/n$ then there exists a non-negative $u_0\in L^1(\R^n)$ such that there is no local integral solution of (\ref{eq:PDE})  that remains in $L^1(\R^n)$ for $t>0$.
\end{lemma}

\Proof
Choose $\alpha<n$ such that $k>(n+2)/ \alpha$. Let $u$ be a local integral solution of (\ref{eq:PDE}) with $f$ constructed as in Section \ref{sec:f} and with initial data $u_0=|x|^{-\alpha}\chi_R$.
Then $u\ge 0$ (since $f\ge 0$ and $u_0\ge 0$) and so  $u(t)\ge S(t)u_0$  for all $t\ge 0$. Hence
\be\label{loweru}
u(t)\ge S(t)u_0+\int_0^tS(t-s)f(S(s)u_0)\,\dee s.
\ee
Integrating over $\R^n$, and noting that $S(t)$ is $L^1$ norm-preserving, we obtain
\bs
\|u(t)\|_{L^1} &=& \int_{\R^n}u(t)\, \dee x \ge\int_{\R^n}S(t)u_0\, \dee x +\int_{\R^n}\int_0^tS(t-s)f(S(s)u_0)\,\dee s\,\dee x\\
&=&\int_{\R^n}u_0\, \dee x +\int_0^t\int_{\R^n}f(S(s)u_0)\, \dee x \, \dee s\\
&\ge &\int_0^t\int_{\R^n} f(S(s)u_0)\, \dee x \, \dee s,
\es
using Fubini's Theorem to change the order of integration.

Now choose $\beta\in(0,1/2)$ such that $k>(n\beta+1)/\alpha\beta$, and write $w(x,t)=[S(t)u_0](x)$. Since $\phi_{i+1}=\phi_i^k$,
\begin{align*}
\int_{\R^n} f(w(x,s))\,\dee x &\ge\int_{\R^n}\ft (w(x,s))\,\dee x\\
&\ge |\{x:\ w(x,s)\ge\phi_i\}|(\phi_{i+1}-\phi_i)\\
&= |\{x:\ w(x,s)\ge\phi_i\}|(\phi_i^k-\phi_i)\\
&\ge |\{x:\ w(x,s)\ge\phi_i\}|\frac{\phi_i^k}{2},\label{eq:int}
\end{align*}
recalling (\ref{eq:phi}).
Now take $i$ sufficiently large so that $(\phi_i/M)^{-1/\alpha\beta}\le t$, recalling (\ref{eq:M}). It now follows from Proposition~\ref{prop:ball} that
\ba
\int_0^t\int_{\R^n}f(w(x,s))\,\dee s\,\dee x&\ge \omega_n &\int_0^{(\phi_i/M)^{-1/\alpha\beta}}\frac{\phi_i^k}{2} s^{n\beta}\,\dee s\nonumber\\
&=& c\phi_i^{k-(n\beta+1)/\alpha\beta}.\label{eq:blowupL1}
\ea
Since $k>(n\beta+1)/\alpha\beta$ the right hand side of (\ref{eq:blowupL1}) tends to $+\infty$ as $i\to\infty$. Hence for this choice of $u_0$ there is no integral solution of (\ref{eq:PDE}) that remains in  $L^1(\R^n)$.\EndProof

We can now prove our main non-existence theorem, in which we obtain blow-up in $L^1_{\rm loc}(\R^n)$, for initial data in $L^q(\R^n)$. The proof shows that this blow-up is local to the origin, i.e.\ does not arise from concentration `at infinity'.

Note that for $q>1$ there is a gap between the values of $k$ for which we can guarantee global existence using the growth condition ($k<1+2q/n$) and those $k$ for which we have non-existence of solutions ($k>q+2q/n$).

\begin{theorem}\label{thm:main}
 Let $q\in [1,\infty )$ and suppose that $f$ is constructed as in (\ref{eq:f}) with $k>q(1+\frac{2}{n})$. Then there exists a $u_0\in L^q(\R^n)$ such that (\ref{eq:PDE}) possesses no local integral solution; indeed, any solution $u(t)$ that satisfies (\ref{eq:mild}) is not in $L^1_{\rm loc}(\R^n)$ for any $t>0$.

\label{thm:nonexist}
\end{theorem}

\Proof
 Fix $\tau \in (0,1)$.  We show that if a local integral solution $u(t)$ exists then it cannot be contained in $L^1(B_n(\rho))$ for any $t\in (0,\tau )$, where $B_n(\rho)$ is the ball of radius $\rho >1$. It follows that for $t\in (0,\tau )$ the solution $u(t)$ cannot be in $L^1_{\rm loc}(\R^n )$, and consequently not in $L^p(\R^n)$ for any $1\le p\le\infty$; therefore $u$ cannot be in $L^1_{\rm loc}((0,\tau)\times\R^n)$, since this would imply that $u(t)\in L^1_{\rm loc}(\R^n)$ for almost every $t\in(0,\tau)$. Since the $f$ we have constructed satisfies $f(s)\ge s$ for all $s\ge\phi_0$, it follows that $f(u)$ cannot be in $L^1_{\rm loc}((0,\tau)\times\R^n)$, and so there is no way to make sense of the integral formula in (\ref{eq:mild}). In other words, there can be no local integral solution, contradicting our initial assumption.

     Since $k>q(1+(2/n))=(n+2)/(n/q)$, there exists an $\alpha<n/q$ such that $k>(n+2)/\alpha$. Choose $\beta\in(0,1/2)$ such that $k>(1+n\beta)/\alpha\beta$. Take $u_0=|x|^{-\alpha}\chi_R\in L^q(\R^n)$ as in Proposition~\ref{prop:ball}. For $t\in (0,\tau )$ fixed, choose $i$ sufficiently large such that $t_i:=(\phi_i/M)^{-1/\alpha\beta}\le t$, where $M$ is as in (\ref{eq:M}). Writing $w(x,t)=[S(t)u_0](x)$, by (\ref{loweru}) we have
\begin{align*}
\int_{B_n(\rho)}u(t)\, \dee x &\ge\int_{B_n(\rho)}\int_0^{t_i}[S(t-s)f(w(\cdot,s))](x)\,\dee s\,\dee x\\
&=C\int_0^{t_i}\int_{B_n(\rho)}(t-s)^{-n/2}\int_{\R^n}\e^{-|x-y|^2/4(t-s)}f(w(y,s))\,\,\dee y\, \dee x \, \dee s\\
&\ge C\int_0^{t_i}(t-s)^{-n/2}\int_{|y|\le s^\beta}\int_{B_n(\rho)}\e^{-|x-y|^2/4(t-s)}\frac{\phi_{i+1}}{2}\, \dee x\,\,\dee y \, \dee s\\
&=\frac{C\phi_{i+1}}{2}\int_0^{t_i}\int_{|y|\le s^\beta}\int_{|z|\le \rho/2\sqrt{t-s}} \e^{-|z-(y/2\sqrt{t-s})|^2}\,\,\dee z\, \dee y \, \dee s\\
&=\frac{C\phi_{i+1}}{2}\int_0^{t_i}\int_{|y|\le s^\beta}\int_{|v|\le (\rho -1)/2\sqrt{t-s}} \e^{-|v|^2}\,\,\dee v\, \dee y \, \dee s,
\end{align*}
where in the last line we have substituted $v=z-(y/2\sqrt{t-s})$ so that the domain of integration is certainly larger than $\{|v|\le (\rho -1)/2\sqrt{t-s}\}$.

Now, since $t_i\le 1$ it follows that
$$
\int_{|v|\le (\rho -1)/2\sqrt{t-s}} \e^{-|v|^2}\,\,\dee v\ge\int_{|v|\le (\rho -1)/2}\e^{-|v|^2}\,\dee v=C^{\prime},
$$
and thus
\begin{align*}
\|u(t)\|_{L^1(B_n(\rho))}&\ge c\phi_{i+1}\int_0^{t_i}\int_{|y|\le s^\beta}\,\dee y\,\dee s\\
&=c\phi_{i+1}\int_0^{t_i} s^{n\beta}\,\dee s\\
&=c\phi_{i+1}t_i^{1+n\beta}=c\phi_i^k\phi_i^{-(1+n\beta)/\alpha\beta}.\qedhere
\end{align*}
Letting $i\to\infty$ yields the result.
\EndProof

Note that if we fix a value of $k$, and hence a nonlinearity $f$, then we can take $q$ sufficiently large that the assumption $k>q(1+2/n)$ of the theorem does not hold. Indeed, we saw earlier that for $k<1+2q/n$ we can guarantee the existence of global solutions for $u_0\in L^q(\R^n)$. So the above result leaves open the question of whether one can find an $f$ for which there are initial conditions in any $L^q(\R^n)$, $1\le q<\infty$, that do not lead to local solutions. However, one can easily construct such an $f$ using the same procedure as above, with $\phi_{i+1}=\e^{\phi_i}$. The blow-up argument is then identical to that of Theorem~\ref{thm:main}, replacing $\phi_i^k$ by $\e^{\phi_i}$ in the final line.

\section{Discussion}
\label{}

We have shown that the Osgood condition
$$
\int_1^\infty\frac{1}{f(s)}\,\dee s=\infty,
$$
necessary and sufficient for the global existence of solutions of the scalar ODE $\dot v=f(v)$, is not sufficient for even local existence of the PDE
$$
u_t=\Delta u+f(u)
$$
when the initial data is unbounded.

Two outstanding questions remain open. First, whether one can obtain similar results from the Dirichlet problem, akin to those
in \cite{CZ} (they treat a problem in which the Osgood-type condition is violated). It appears that significantly more refined arguments would be necessary to extend our results to this case, but we would expect such results to hold given the local nature of the blow-up in Theorem~\ref{thm:nonexist}.

More technically, our proof involved the construction of a particular nonlinearity $f$, whose growth was parametrised by $k\in(1,\infty)$. For initial data in $L^1(\R^n)$ we have shown that there is a sharp transition between global existence ($k<1+2/n$) and local non-existence ($k>1+2/n$), although we have not addressed the transition point $k=1+2/n$. For initial data in $L^q(\R^n )$ there is a gap between the range in which we can prove global existence ($k<1+2q/n$) and in which we can show local non-existence ($k>q(1+2/n)$). It would be interesting to understand the behaviour of this model for $1+2q/n\le k\le q(1+2/n)$.


\vspace{1cm}
\noindent{\bf Acknowledgements}
The authors acknowledge support under the following grants:  EPSRC Platform EP/I019138/1 (RL);  EPSRC Leadership Fellowship  EP/G007470/1 (JCR);  EPSRC Doctoral Prize EP/P50578X/1 (MS).





%

\end{document}